\documentclass[12pt]{article}
\newtheorem{theorem}{Theorem}
\begin{document}
\title{The Redheffer matrix of a partially ordered set}
\author{Herbert S. Wilf\\University of Pennsylvania\\Philadelphia, PA 19104-6395}
\maketitle
\begin{abstract}
R. Redheffer described an $n\times n$ matrix of 0's and 1's the size of whose determinant is connected to the Riemann Hypothesis. We describe the permutations that contribute to its determinant and evaluate its permanent in terms of integer factorizations. We generalize the Redheffer matrix to finite posets that have a 0 element and find the analogous results in the more general situation.
\end{abstract}
\section{Introduction}
In 1977, R. Redheffer described a matrix that is closely connected to the Riemann Hypothesis (RH). Let $R_n$ be the $n\times n$ matrix whose $(i,j)$ entry is 1 if $i\backslash j$ or if $j=1$, and otherwise is 0, for $1\le i,j\le n$. For example,
\[R_8=\left(
\begin{array}{cccccccc}
1&1&1&1&1&1&1&1\\1&1&0&1&0&1&0&1\\1&0&1&0&0&1&0&0\\1&0&0&1&0&0&0&1\\1&0&0&0&1&0&0&0\\1&0&0&0&0&1&0&0\\1&0&0&0&0&0&1&0\\1&0&0&0&0&0&0&1
\end{array}\right)
\]
He showed that the proposition ``for every $\epsilon>0$ we have $|\det{R_n}|=O(n^{\frac12+\epsilon})$'' is equivalent to RH. More precisely, he showed that
\begin{equation}
\label{eq:mert}
\det{(R_n)}=\sum_{k=1}^n\mu(k),
\end{equation}
where $\mu$ is the classical M\"obius function, and the equivalence of the $O(n^{\frac12+\epsilon})$ growth bound of the right side of (\ref{eq:mert}) to the RH is well known.

Here we will first describe the permutations of $n$ letters that contribute to the determinant of $R_n$, i.e., the permutations that do not hit any 0 entries in the matrix. Then we will count those permutations, which is to say that we will evaluate the permanent of the Redheffer matrix. It turns out that this permanent is also nicely expressible in terms of well known number theoretic functions.

After that we will generalize the Redheffer matrix to posets other than the positive integers under divisibility,   and find an application in the case of the Boolean lattice.
\section{The permanent}
Which permutations of $n$ letters contribute a $\pm 1$ to the determinant above? Let
\[\tau=R_{i_1,1}R_{i_2,2}\dots R_{i_n,N}\]
 be a nonvanishing term of that determinant. Fix some integer $j_1$, $2\le j_1\le n$. Then in the term $\tau$ there is a factor $R_{j_1,j_2}$ for some unique $j_2$. If $j_2\neq j_1$, then there is also a factor $R_{j_2,j_3}$ for some $j_3$, etc. Finally, we will have identified a collection of nonvanishing factors $\sigma=R_{j_1,j_2}R_{j_2,j_3}\dots R_{j_k,j_1}$ in the term $\tau$. By the definition of the matrix, we must have $j_1<j_2\le j_3\le \dots \le j_k\le j_1$, which is a contradiction because $j_1$ was assumed to be $>1$.

It follows that in any collection $\sigma$ of contributing factors, we either have $j_2=j_1$, i.e., $\sigma$ has just a single factor in it, namely a diagonal element of the matrix, or else $j_1=1$. Suppose $j_1=1$. Then the collection $\sigma$ is of the form $R_{1,j_2}R_{j_2,j_3}\dots R_{j_k,1}$, and for this to give a nonzero contribution what we need is that $j_2\backslash j_3\backslash \dots \backslash j_k$, i.e., the sequence $j_2,j_3,\dots ,j_k$ forms a chain under divisibility.

We can therefore match nonvanishing contributions to the determinant of $R_n$ with permutations in which the cycle that contains 1 is a division chain, and the other cycles are all fixed points.

To phrase this in more a traditional number theoretic way, recall that an \textit{ordered factorization} of an integer $m$ is a representation $m=a_1a_2\dots a_k$, in which all $a_i$s are $\ge 2$ and the order of the factors is important. Now in our case, the successive quotients
\[j_2/j_1,j_3/j_2,\dots j_k/j_{k-1}\qquad (j_1=1)\]
 are an ordered factorization of some integer (namely $j_k$) which is $\le n$. The number of contributing permutations is therefore equal to the number of all ordered factorizations of all positive integers $\le n$, plus 1 more to account for the empty factorization. Hence we have the following result.
\begin{theorem}
The permanent of the Redheffer matrix $R_n$ is $1+\sum_{k=1}^nf(k)$, where $f(k)$ is the number of ordered factorizations of the integer $k$. The permutations that contribute to this permanent are those in which there is just one cycle of length $>1$, the letter $1$ lives in that cycle, and the elements of that cycle form a division chain.
\end{theorem}
The values of these permanents form sequence A025523 in Sloane's database \cite{ns} of integer sequences. For $n=1$ to $10$ their values are $1,2,3,5,6,9,10,14,16,19$. It is known \cite{hh} that this sequence grows like $Cn^a$ for large $n$, where $a=\zeta^{-1}(2)=1.73..$.

It is now easy to give another proof of Redheffer's evaluation of the determinant. If the cycle that contains 1 contains $k$ letters altogether, the highest of which is $r$, then the contributing permutation has $n-k+1$ cycles on $n$ letters, so its sign is $(-1)^{k-1}$. The contribution of all such permutations in which the highest letter of the cycle that contains 1 is $r$ is $\sum_{\phi}(-1)^{k(\phi)}$, extended over all ordered factorizations $\phi$ of $r$, where $k(\phi)$ is the number of factors in $\phi$. It is a known result from number theory (for a bijective proof see \cite{gkw}) that this sum is
\[\sum_{\phi}(-1)^{k(\phi)}=\mu(r),\]
from which the evaluation (\ref{eq:mert}) follows by summing on $r$.
\section{Generalizations}
The analysis of the preceding section can be carried out in general posets. Let $(\cal{S},\preceq)$ be a finite poset that has a 0 element, and suppose the elements of ${\cal S}$ have been labeled by the positive integers so that the $\zeta$ matrix of ${\cal S}$ is upper-triangular. We define the Redheffer matrix $R(S)$ of $S$ to be the result of replacing the first column (i.e., the column that is labeled by the 0-element) of the $\zeta$-matrix of $S$ by a column of all 1's.

By the argument of the preceding section, the permanent of $R(S)$ is the number of $\preceq$-chains in $S$ that contain the $0$ element. The permutations of $[\,|S|\,]$ that contribute to the permanent are those all of whose cycles are fixed points except for the cycle that contains 0, which must be a chain in $S$. The determinant of $R(S)$ is the sum of $(-1)^{L(C)}$ over all chains $C$ in ${\cal P}-\{0\}$, where $L(C)$ is the length of the chain $C$. If we group the terms of this sum according to the largest element of each chain $C$, then the contributions whose largest element is some fixed $x$ sum up to $\mu(0,x)$, where $\mu$ is the M\"obius function of the poset. If we sum over $x$ we find that the determinant of the general Redheffer matrix is $\sum_x\mu(0,x)$.
\begin{theorem}
Let $R=R(S)$ be the Redheffer matrix of a finite poset $S$ that contains a 0 element. Then the permanent of $R$ is the number of chains of $S$ that contain the 0 element, and the determinant of $R$ is $\sum_{x\in S}\mu(0,x)$, where $\mu$ is the M\"obius function of $S$.
\end{theorem}

 If the poset has a ``1'' element then this sum is 0, and the Redheffer matrix is singular. Thus in the Boolean lattice ${\cal B}_n$ on $n$ elements, for example, the $2^n\times 2^n$ Redheffer matrix is singular. Its permanent is the number of chains in ${\cal B}_n$ that contain $\{\emptyset\}$. These numbers
\[1,2,6,26,150,1082,\dots\]
form sequence number A00629 in Sloane's database.

It is easy to write out the inverse of the generalized Redheffer matrix. This follows at once from the formula for finding the inverse of a matrix that differs from one of known inverse by a matrix of rank 1. The formula is
\[(B+\mathbf{u}\mathbf{v}^T)^{-1}=B^{-1}-\frac{B^{-1}\mathbf{u}\mathbf{v}^TB^{-1}}{1+(\mathbf{v},B^{-1}\mathbf{u})}.\]
It yields, in our case,
\[(R^{-1})_{x,y}=\mu(x,y)-\frac{\mu(0,y)\sum_{z\neq 0}\mu(x,z)}{\sum_z\mu(0,z)}.\qquad(x,y\in S)\]


\newpage

\begin{thebibliography}{aaa}
\bibitem{bfp} W. W. Barrett, R. W. Forcade and A. D. Pollington, On the spectral radius of a (0,1) matrix related to Mertens' function, Linear Algebra Appl. \textbf{107} (1988), 151–-159;
\bibitem{ba} Wayne W. Barrett and Tyler J. Jarvis, Spectral properties of a matrix of Redheffer, \textit{Directions in matrix theory} (Auburn, AL, 1990), Linear Algebra Appl. \textbf{162/164} (1992), 673--683.
\bibitem{gkw} Richard Garfield, Donald E. Knuth, and Herbert S. Wilf, A bijection for ordered factorizations, J. Combin. Theory Ser. A \textbf{54} (1990), no. 2, 317--318.
\bibitem{hu} Stephen P. Humphries, Cogrowth of groups and a matrix of Redheffer, Linear Algebra Appl. \textbf{265} (1997), 101--117.
\bibitem{hh} Hsien-Kuei Hwang, Distribution of the number of factors in random ordered factorizations of integers, J. Number Theory \textbf{81} (2000), no. 1, 61--92.
\bibitem{ja} Tyler J. Jarvis, A dominant negative eigenvalue of a matrix of Redheffer, Linear Algebra Appl. \textbf{142} (1990), 141--152.
\bibitem{re} R. M. Redheffer, Eine explizit l\"osbare Optimierungsaufgabe, \textit{Internat. Schriftenreihe Numer. Math.} \textbf{36} (1977).
\bibitem{ns} Neil J. A. Sloane, The On-Line Encyclopedia of Integer Sequences, on the web at \texttt{<http://www.research.att.com/~njas/sequences>}.
\bibitem{rv2} R. C. Vaughan, On the eigenvalues of Redheffer's matrix, I, \textit{Number theory with an emphasis on the Markoff spectrum} (Provo, UT, 1991), 283--296, Lecture Notes in Pure and Appl. Math., \textbf{147}, Dekker, New York, 1993.
\bibitem{rv1} R. C. Vaughan, On the eigenvalues of Redheffer's matrix, II, J. Austral. Math. Soc. Ser. A \textbf{60} (1996), no. 2, 260--273.

\end{thebibliography}
\end{document}